\documentclass[11pt]{amsart}

\newtheorem{theorem}{Theorem}[section]

\newtheorem{lemma}{Lemma}[section]

\theoremstyle{definition}
\newtheorem{remark}{Remark}[section]

\numberwithin{equation}{section}

\newcommand{\hilbert}{{\mathcal H}}
\newcommand{\unitop}{{\mathfrak W}}
\newcommand{\costr}{{\mathfrak C}}

\begin{document}

\title[Zeta functional equation]
{A general statement of the functional 
equation for the Riemann zeta-function}

\author{Luis B\'{a}ez-Duarte}

\date{21 May 2002}
\maketitle

\section{Introduction}
This is a completely formal exposition. All lemmas and theorems are valid for
sufficiently large, dense subspaces of  $\hilbert:=L_2((0,\infty),dx)$.  We borrow
\textit{in toto} the notations from my paper on invariant unitary operators
\cite{iuo}. We denote the Mellin transform of $f$ by

$$
f^{\wedge}(s):=\int_0^\infty x^{s-1}f(x)dx,
$$
which is at least defined on the critical line as a Cauchy principal value in the
$\hilbert$-sense. 
\ \\

Following Burnol \cite{burnol} denote the M\"{u}ntz-modified Poisson summation
operator by $P$, which acts formally on the set of all complex
functions $f$ on $(0,\infty)$ according to

$$
Pf(x):=\sum_{n=1}^\infty f(nx)-\frac{1}{x}\int_0^\infty f(t)dt.
$$
\\ 
The following lemma is obvious:

\begin{lemma}\label{l1}
For any invariant operator $W$ we have $PW=WP$.
\end{lemma}
\ \\
We take for granted the following facts which apply
symultaneously to any
\textit{good} function $f\in\hilbert$.

\begin{lemma}[M\"{u}ntz's formula]\label{l2}
$$
(Pf)^{\wedge}(s)=\zeta(s)f^{\wedge}(s),
$$
\end{lemma}

\begin{lemma}\label{l3}
$$
(Sf)^{\wedge}(s)=f^{\wedge}(1-s).
$$
\end{lemma}

\begin{lemma}\label{l4}
$$
(\costr f)^{\wedge}(s)=\phi(s)f^{\wedge}(1-s),
$$
where
$$
\phi(s):=2^{1-s}\pi^{-s}\cos\frac{\pi s}{2}\Gamma(s).
$$
\end{lemma}

\section{General statements of the functional equation}
We recall Theorem 2.1 of \cite{iuo}:

\begin{theorem}\label{exist}
For any skew root $T$, and any real-valued generator $g$ there is a unique, bounded,
invariant operator $W = \unitop (T,g)$ with $Wg = Tg$. W is unitary and satisfies
$(TW)^{2}=(WT)^{2}=I$, $W^{-1}=TWT$. 
\end{theorem}

Why call the following triviality a ``theorem" is a question left to the reader. Do
note that the operator $W$ depends explicitly on $f$.

\begin{theorem}[Main identity]\label{main}
Let $f$ be good and suppose that $Pf$ is a real valued generator. Let $T$ be a skew
root and $W=\unitop(T,Pf)$. Then
$$
TPf=PWf.          
$$
\end{theorem}

\begin{proof}
Use Theorem \ref{exist} and the commutativity Lemma \ref{l1}.
\end{proof}

We specialize the main identity to give two general statements of the functional
equation of the Riemann zeta function.

\begin{theorem}[Functional equation I]\label{fe1}
Let $f$ be good and suppose that $Pf$ is a real valued generator. Let 
$W=\unitop(S,Pf)$ and $h=Wf$. Then

\begin{equation}\label{funceq1}
\zeta(1-s)f^{\wedge}(1-s)=\zeta(s)h^{\wedge}(s).          
\end{equation}
\end{theorem}

\begin{proof}
The main identity is now $SPf=Ph$, where we apply the Mellin transform using Lemmas
\ref{l2}, \ref{l3}.
\end{proof}

\begin{theorem}[Functional equation II]\label{fe2}
Let $f$ be good and suppose that $Pf$ is a real valued generator. Let 
$W=\unitop(\costr,Pf)$ and $h=Wf$. Then

\begin{equation}\label{funceq2}
\zeta(1-s)\phi(s)f^{\wedge}(1-s)=\zeta(s)h^{\wedge}(s).          
\end{equation}
\end{theorem}

\begin{proof}
The main identity is now $\costr Pf=Ph$, where we apply the Mellin transform using
Lemmas
\ref{l2}, \ref{l4}.
\end{proof}

For an application of Theorem \ref{fe1} choose $f=\chi$, then $Pf(x)=-\rho(1/x)$ and
$W=U$, so that

$$
h(x)=\frac{\sin(2\pi x)}{\pi x}.
$$
Therefore

$$
f^{\wedge}(s)=\frac{1}{s}, \ \ \ h^{\wedge}(s)=-2^{1-s}\pi^{-s} \cos\frac{\pi
s}{2}\Gamma(s-1),
$$
so (\ref{funceq1}) expresses the functional equation as

$$
\frac{\zeta(1-s)}{1-s}=-2^{1-s}\pi^s \cos\frac{\pi s}{2}\Gamma(s-1)\zeta(s).
$$
\\

I don't yet have a genuine application  of Theorem \ref{fe2}. However let as usual
$M(x):=\sum_{n\leq x}\mu(n),$ and choose 

$$
f(x):=M(1/x).
$$
\\
I have shown in \cite{lbd} that this $f\not\in\hilbert$; on the other hand note that
$Pf=\chi$. So again $W=\unitop(\costr,Pf)=\unitop(\costr,\chi)=U$. If one
blithely insists on applying the theorem one obtains an intriguing equation from
(\ref{funceq2}).

\begin{remark}
Quite frankly I am chagrined that I have not been included other genuine
examples. More important perhaps is this \textbf{question}: how does this relate to
Burnol's view of the functional equation in \cite{burnol}?
\end{remark}

\bibliographystyle{amsplain}

\ \\
\noindent Luis B\'{a}ez-Duarte\\
Departamento de Matem\'{a}ticas\\
Instituto Venezolano de Investigaciones Cient\'{\i}ficas\\
Apartado 21827, Caracas 1020-A\\
Venezuela\\
\email{lbaezd@cantv.net}

\end{document}